\documentclass[11pt]{article}

\usepackage[margin=1in]{geometry}
\usepackage{amsmath,amssymb,amsthm,mathtools}
\usepackage{booktabs}
\usepackage{longtable}
\usepackage[hidelinks]{hyperref}
\usepackage[nocompress]{cite}

\newtheorem{theorem}{Theorem}[section]
\newtheorem{lemma}[theorem]{Lemma}
\newtheorem{proposition}[theorem]{Proposition}
\newtheorem{corollary}[theorem]{Corollary}
\newtheorem{conjecture}[theorem]{Conjecture}
\theoremstyle{definition}

\theoremstyle{remark}
\newtheorem{remark}[theorem]{Remark}

\newcommand{\Ztwo}{\mathbb Z_2}
\newcommand{\Qtwo}{\mathbb Q_2}

\title{A closed-form nested-sum formula for fixed-size matchings on square tori}
\author{Yong Kong\\
Department of Biostatistics, Yale School of Public Health\\
New Haven, Connecticut 06520, USA\\
\texttt{yong.kong@yale.edu}}
\date{}

\begin{document}
\maketitle

\begin{abstract}
Let $a(N,s)$ denote the number of matchings of size $s$ on the square
torus $C_m\mathbin{\square}C_n$, where $N=mn$ and $m,n>s+1$.
In statistical physics this enumeration belongs to the monomer--dimer
problem; in computational complexity, its general planar-graph form is
\#P-complete.
Previous work proved that $a(N,s)$ is a polynomial of degree $s$ in
$N$.  We derive a finite closed-form nested-sum formula for the entire
polynomial.  The construction starts from the per-vertex connected
coefficients $b_k$ in the Mayer expansion and applies an iterated
reversion--logarithm transformation obtained from Lagrange inversion,
producing rational constants $g_i$.  Writing $v_2$ for the 2-adic
valuation, we prove the divisibility property $k b_k\in2\mathbb Z$ and
use it to establish the exact valuation
\[
  v_2(g_i)=1-2^{i-2}\qquad(i\geq2)
\]
for the rational constants generated by the transformation.  In
particular, every stage is nondegenerate, so the nested formula is
defined at arbitrary depth.  The valuation theorem gives a doubly
exponential lower bound on the reduced denominators of $g_i$.
Extensive exact computations further provide strong evidence that the
magnitudes $|g_i|$ themselves grow doubly exponentially.  The formula
contains exactly $2^{s-1}$ chain-indexed summands and provides a
concrete starting point for further arithmetic and complexity analysis
of fixed-dimer enumeration.  The underlying transformation is
model-independent and also supplies nested formulas for other
finite-range hard-core models, including aligned $k$-mers, whenever the
corresponding transformation sequence is nondegenerate.
\end{abstract}

\noindent\textbf{Keywords.}
Monomer--dimer model; matching polynomial; exact enumeration; Mayer
cluster expansion; Lagrange inversion; formal power series; 2-adic
valuation; double-exponential growth.

\medskip
\noindent\textbf{2020 Mathematics Subject Classification.}
05A15; 05C70; 82B20; 68Q25.

\section{Introduction}\label{sec:introduction}

A dimer on a graph is an edge, and a configuration of $s$ nonoverlapping
dimers is a matching of size $s$.  In equilibrium statistical mechanics,
the same objects form the monomer--dimer model: vertices incident with a
selected edge are covered by dimers and the remaining vertices are
monomers.  The model is one of the classical hard-core lattice systems
and has been studied from statistical physics, combinatorics, and
computational complexity
\cite{Kasteleyn1961,TemperleyFisher1961,Fisher1961,HeilmannLieb1972,
  Jerrum1987}.
Close-packed dimers on planar graphs admit Pfaffian methods, whereas the
general monomer--dimer enumeration problem is computationally much more
difficult: counting matchings in a general planar graph is
\#P-complete \cite{Kasteleyn1961,TemperleyFisher1961,Jerrum1987}.

We consider the graph
\[
  G_{m,n}=C_m\mathbin{\square}C_n,
\]
the Cartesian product of two cycles.  It is the $m\times n$ square
lattice with periodic boundary conditions in both directions.  The
number of its vertices is denoted by
\[
  N=mn.
\]
For a nonnegative integer $s$, let $a_{m,n}(s)$ be the number of
matchings of size $s$ in $G_{m,n}$.  A recent general polynomiality
theorem for finite-range hard-core models implies that, whenever
$m,n>s+1$, this number depends on $m$ and $n$ only through $N$ and is a
polynomial of degree $s$ in $N$ \cite{Kong2026}.  In the stable range we
therefore write the count as $a(N,s)$.

The polynomiality theorem does not by itself expose the coefficients of
this polynomial.  The purpose of the present paper is to determine their
complete nested structure.  Our starting point is the familiar
linked-cluster principle: the logarithm of a partition function retains
only connected contributions.  In statistical mechanics this is the
Mayer cluster expansion \cite[Chapter~5]{Hill1956}; see also
\cite{MayerMayer1940,Ruelle1969}.  We give a short finite-variable proof
adapted to the matching polynomial because it supplies the precise
stability and arithmetic information needed below.  Exact and
series-based approaches to the square-lattice problem have a long
history \cite{Gaunt1969,Kong2006,Kong2006b,Kong2006c,Kong2007,Kong2024}.

The main contributions are as follows.
\begin{enumerate}
  \item We obtain a finite closed-form expression for $a(N,s)$ indexed
  by all strictly decreasing chains beginning at $s$.  It contains
  $2^{s-1}-1$ nonempty chain terms in addition to the leading term and
  requires neither interpolation in $N$ nor a recurrence for the
  unknown count.
  \item We identify and compute the rational sequence $g_i$ governing
  every level of the nesting.  Exact calculations through $i=26$ show
  that $|g_i|^{1/2^i}$ stabilizes near $1.3197289$, giving striking
  evidence that $|g_i|$ grows doubly exponentially.  This Archimedean
  growth is stated as a conjecture rather than used as an assumption.
  \item As a supporting arithmetic result, we define stable connected
  coefficients $b_k$ and prove
  \[
    c_k:=k b_k\in2\mathbb Z\qquad(k\geq1).
  \]
  \item We introduce the iterated reversion--logarithm transformation and
  derive its exact coefficient action by formal Lagrange inversion.
  \item We prove that every $g_i$ is nonzero and establish the exact
  2-adic valuation
  \[
    v_2(g_i)=1-2^{i-2}\qquad(i\geq2).
  \]
  Thus the nesting is defined at arbitrary depth, and the reduced
  denominators of $g_i$ contain an exactly determined doubly
  exponential power of $2$.
\end{enumerate}

The first two items are the principal structural findings; the final
three provide the transformation and nonvanishing mechanism that make
the closed formula rigorous.  A preliminary version of this work was
posted on arXiv.\footnote{Y.~Kong, ``Intractable enumeration problems
are like Russian nesting dolls: structural properties of monomer--dimer
coverings on two-dimensional quadratic lattices,'' arXiv:2609.19231
(2026), \url{https://arxiv.org/abs/2609.19231}.}  The present manuscript
focuses the exposition on the closed formula and supplies a standalone
proof of its nondegeneracy.

The word \emph{closed-form} is used here in the following precise sense.
For fixed $s$, the answer is a finite sum over explicitly specified
integer chains; all weights are finite rational expressions in the
universal constants $g_1,\ldots,g_s$, and these constants are determined
by explicit coefficient extraction from the stable Mayer series.  The
formula requires neither interpolation in $N$ nor a recurrence involving
the unknown count itself.

Section~\ref{sec:stable} establishes the stable Mayer coefficients and
their divisibility.  Section~\ref{sec:transform} develops the formal
transformation.  Section~\ref{sec:valuation} proves nonvanishing at every
depth.  The closed nested formula is proved in
Section~\ref{sec:closed}.  Section~\ref{sec:numerics} records exact data
and the observed growth of $g_i$, and Section~\ref{sec:complexity}
discusses the scope and limitations of the resulting complexity
interpretation.

\section{Stable Mayer coefficients}\label{sec:stable}

Define the ordinary matching generating polynomial
\begin{equation}\label{eq:finiteZ}
  Z_{m,n}(z)=\sum_{q\geq0}a_{m,n}(q)z^q.
\end{equation}
Because $G_{m,n}$ is finite, the sum in \eqref{eq:finiteZ} is a
polynomial and $Z_{m,n}(0)=1$.  Consequently its formal logarithm belongs
to $z\mathbb Q[[z]]$.

The stable polynomiality result used in this paper is the following
specialization of \cite[Lemma~1 and Theorem~1]{Kong2026}.

\begin{proposition}[Stable polynomiality]\label{prop:polynomiality}
For every $s\geq0$ and $m,n>s+1$, there is a polynomial $P_s(X)\in
\mathbb Q[X]$ of degree $s$, independent of $m$ and $n$, such that
\[
  a_{m,n}(s)=P_s(mn).
\]
\end{proposition}

We next refine Proposition~\ref{prop:polynomiality} at the level of the
logarithm of the partition function.

\begin{lemma}[Stable connected coefficients]\label{lem:stable}
There are rational numbers $b_k$, independent of $m$ and $n$, such that
for every $k\geq1$ and every $m,n>k+1$,
\begin{equation}\label{eq:stablebk}
  [z^k]\log Z_{m,n}(z)=N b_k,
  \qquad N=mn.
\end{equation}
Moreover,
\begin{equation}\label{eq:firstb}
  b_1=2,\qquad b_2=-7,
\end{equation}
and
\begin{equation}\label{eq:divisibility}
  c_k:=k b_k\in2\mathbb Z\qquad(k\geq1).
\end{equation}
\end{lemma}

\begin{proof}
Assign a separate commuting variable $x_e$ to every edge $e$ of
$G_{m,n}$ and form the multivariate matching polynomial
\[
  Z_{m,n}(\boldsymbol{x})=
  \sum_{M}\prod_{e\in M}x_e,
\]
where the sum is over all matchings $M$, including the empty matching.
Suppose that the distinct edges in the support of a monomial split into
two nonempty collections with no common endpoint between the
collections.  After every variable outside this support is set to zero,
the matching polynomial factors into the product of the matching
polynomials of the two collections.  Its logarithm is therefore a sum,
and the coefficient of the original mixed monomial in
$\log Z_{m,n}(\boldsymbol{x})$ is zero.  Thus every monomial of total
degree $k$ contributing to the logarithm has connected support containing
at most $k$ distinct lattice edges.  This is the finite formal version of
the Mayer connected-cluster principle.

If $m,n>k+1$, such a support cannot wind around either periodic
direction.  It lifts uniquely, after a root edge and a translation have
been fixed, to a local configuration in the infinite square lattice.
The coefficient of the monomial and the list of rooted local shapes are
independent of the periods.  After choosing whether a distinguished edge
is horizontal or vertical, there are only finitely many lifted local
types.  A nonempty, connected, nonwinding support has trivial stabilizer
under torus translations: a nontrivial translation fixing it would
repeat the support around a periodic direction, requiring at least the
corresponding period in its connected edge support.  Hence every such
oriented type has exactly $N$ translates.  Its coefficient in the
formal logarithm depends only on the lifted local type.  Summing over
the finite list of types and then setting every $x_e=z$ therefore gives
\eqref{eq:stablebk}, with the remaining finite rational sum denoted by
$b_k$.

There are $2N$ edges.  In the range used here, $m,n>3$, every edge has
six distinct neighboring edges that share an endpoint with it.  Each
conflicting unordered pair is counted from both of its edges, so there
are $(2N)6/2=6N$ such pairs.  Hence
\[
 a(N,1)=2N,
 \qquad
 a(N,2)=\binom{2N}{2}-6N=2N^2-7N.
\]
Comparison with the logarithm of
$1+a(N,1)z+a(N,2)z^2+O(z^3)$ yields \eqref{eq:firstb}.

It remains to prove \eqref{eq:divisibility}.  Fix $k$ and choose a square
torus $G_{L,L}$ with $L>k+1$.  Write $Z_G(z)$ for its matching
polynomial.  For an edge $e=\{u,v\}$, let $G-\{u,v\}$ denote the induced
graph obtained by deleting its endpoints.  Counting a matching together
with one distinguished edge gives the formal identity
\begin{equation}\label{eq:pointed}
  \frac{zZ_G'(z)}{Z_G(z)}
  =\sum_{e=\{u,v\}\in E(G)}
    z\frac{Z_{G-\{u,v\}}(z)}{Z_G(z)}.
\end{equation}
Translations and a quarter-turn act transitively on the $2N$ edges of
the square torus.  Every quotient on the right of
\eqref{eq:pointed} is the same series $Q(z)$.  Since both polynomials in
the quotient have integer coefficients and $Z_G(0)=1$, we have
$Q(z)\in\mathbb Z[[z]]$.  The right side of \eqref{eq:pointed} therefore
belongs to $2N z\mathbb Z[[z]]$.  On the other hand, by
\eqref{eq:stablebk}, the coefficient of $z^k$ on the left is $kNb_k$.
It follows that $kb_k\in2\mathbb Z$.
\end{proof}

Define the stable Mayer series
\begin{equation}\label{eq:Bdef}
  B(z)=\sum_{k\geq1}b_kz^k
      =\sum_{k\geq1}\frac{c_k}{k}z^k.
\end{equation}
For a formal variable $X$, define
\begin{equation}\label{eq:Pdef}
  P_s(X)=[z^s]\exp\bigl(XB(z)\bigr).
\end{equation}
If $m,n>s+1$, then \eqref{eq:stablebk} through degree $s$ implies
\begin{equation}\label{eq:stableP}
  a(N,s)=P_s(N).
\end{equation}
Only $b_1,\ldots,b_s$ enter \eqref{eq:Pdef}, so this is a finite formal
identity for every fixed $s$.  It also shows directly that $P_s$ has
degree $s$, zero constant term for $s\geq1$, and leading coefficient
$2^s/s!$.

\section{The reversion--logarithm transformation}\label{sec:transform}

All operations in this section take place in formal power-series rings
over $\mathbb Q$.  No analytic convergence is required.  For a formal
series
\begin{equation}\label{eq:Adef}
  A(z)=gz+a_2z^2+a_3z^3+\cdots,
  \qquad g\ne0,
\end{equation}
define
\begin{equation}\label{eq:Fdef}
  F_A(z)=\frac{A(z)}{gz}.
\end{equation}
Thus $F_A(0)=1$.  Let $T_A(w)\in w\mathbb Q[[w]]$ be the unique
compositional inverse of
\begin{equation}\label{eq:WA}
  W_A(z)=zF_A(z);
\end{equation}
that is, $W_A(T_A(w))=w$ and $T_A(W_A(z))=z$.

For a sign $\epsilon\in\{1,-1\}$, define the
\emph{reversion--logarithm transform}
\begin{equation}\label{eq:transform}
  \mathcal L_\epsilon(A)(t)
  =\log F_A\bigl(T_A(\epsilon gt)\bigr).
\end{equation}
Whenever a single transformation is under discussion, we use the
explicit notation
\begin{equation}\label{eq:Aplusdef}
  A^+(t):=\mathcal L_\epsilon(A)(t)
          =g^+t+a_2^+t^2+a_3^+t^3+\cdots.
\end{equation}
Thus $A^+$ denotes the transformed series, $g^+$ denotes its linear
coefficient, and $a_r^+$ denotes its coefficient of $t^r$.

Starting from $A_1=B$, define recursively
\begin{equation}\label{eq:iteration}
  A_{i+1}=\mathcal L_{\epsilon_i}(A_i),
  \qquad
  \epsilon_1=-1,
  \qquad
  \epsilon_i=1\quad(i\geq2),
\end{equation}
and write
\begin{equation}\label{eq:stagecoeff}
  A_i(z)=\sum_{r\geq1}a_{r,i}z^r,
  \qquad
  g_i:=a_{1,i}=A_i'(0).
\end{equation}
At this point the iteration is defined as long as every $g_i$ is
nonzero.  Section~\ref{sec:valuation} proves that it is defined for all
$i$.  From $b_1=2$ and $b_2=-7$ we already obtain
\begin{equation}\label{eq:firstg}
  g_1=2,
  \qquad
  g_2=\epsilon_1b_2=7.
\end{equation}

The next lemma is the formal mechanism behind the nesting.  The formal
Lagrange--B\"urmann identity used below is standard
\cite[Chapter~I.5]{FlajoletSedgewick2009}.

\begin{lemma}[Coefficient transformation]\label{lem:lagrange}
Let $A$ be as in \eqref{eq:Adef}, let $\epsilon\in\{1,-1\}$, and define
$A^+=\mathcal L_\epsilon(A)$ by \eqref{eq:Aplusdef}.  For integers
$n>j\geq1$,
\begin{align}
 \label{eq:coefffirst}
 [X^{n-j}][z^n]e^{XA(z)}
 &=\frac{\epsilon^jg^{n-2j}}{(n-j-1)!}
   \frac{[t^j]e^{nA^+(t)}}{n}\\
 \label{eq:coeffsecond}
 &=\frac{\epsilon^jg^{n-2j}}{(n-j-1)!}
   \sum_{\ell=0}^{j-1}
   \frac{[t^j]A^+(t)^{j-\ell}}{(j-\ell)!}
   n^{j-\ell-1}.
\end{align}
If $g^+=(A^+)'(0)\ne0$, the sum in
\eqref{eq:coeffsecond} is a polynomial in $n$ of degree exactly $j-1$,
with leading coefficient $(g^+)^j/j!$.
\end{lemma}

\begin{proof}
Taking the coefficient of $X^{n-j}$ first gives
\begin{equation}\label{eq:coefficientstart}
 [X^{n-j}][z^n]e^{XA(z)}
 =\frac{g^{n-j}}{(n-j)!}[z^j]F_A(z)^{n-j}.
\end{equation}
Set $w=W_A(z)=zF_A(z)$, so that $z=T_A(w)$.  Apply the formal
Lagrange--B\"urmann formula to $\Phi(z)=F_A(z)^n$:
\begin{align*}
 [w^j]F_A(T_A(w))^n
 &=\frac{n}{j}[z^{j-1}]F_A'(z)F_A(z)^{n-j-1}\\
 &=\frac{n}{n-j}[z^j]F_A(z)^{n-j}.
\end{align*}
Rearranging gives
\begin{equation}\label{eq:residueidentity}
  \frac{[z^j]F_A(z)^{n-j}}{n-j}
  =\frac{[w^j]F_A(T_A(w))^n}{n}.
\end{equation}

By \eqref{eq:transform},
\[
  F_A(T_A(\epsilon gt))^n=e^{nA^+(t)}.
\]
Substituting $w=\epsilon gt$ into \eqref{eq:residueidentity} and then
into \eqref{eq:coefficientstart} proves \eqref{eq:coefffirst}.
Expanding $e^{nA^+(t)}$ proves \eqref{eq:coeffsecond}.  Because
$A^+(0)=0$, only the powers $1,\ldots,j$ contribute to $[t^j]$.
The contribution from the $j$th power is
$(g^+)^jn^{j-1}/j!$, proving the last assertion.
\end{proof}

The following coefficient form of the transformation will be used in
the valuation argument.

\begin{lemma}[Exact transformed coefficients]\label{lem:transformcoeff}
Let $A$, $A^+$, and $\epsilon$ be as in
Lemma~\ref{lem:lagrange}, and assume $a_2\ne0$.  Define the normalized
variable
\begin{equation}\label{eq:udef}
  u=\frac{a_2}{g}z
\end{equation}
and the normalized series
\begin{equation}\label{eq:Fhatdef}
  \widehat F(u)
  :=F_A\!\left(\frac{g}{a_2}u\right).
\end{equation}
Then, for every $r\geq1$,
\begin{equation}\label{eq:valuationtransform}
  a_r^+=[t^r]A^+(t)
  =-\frac{(\epsilon a_2)^r}{r}
    [u^r]\widehat F(u)^{-r}.
\end{equation}
In particular, $g^+=\epsilon a_2$.
\end{lemma}

\begin{proof}
Apply the Lagrange--B\"urmann formula to
$\log F_A(T_A(w))$.  Since $T_A(w)$ is the inverse of
$w=zF_A(z)$,
\[
 [w^r]\log F_A(T_A(w))
 =\frac1r[z^{r-1}]
   \frac{F_A'(z)}{F_A(z)}F_A(z)^{-r}
 =-\frac1r[z^r]F_A(z)^{-r}.
\]
Replacing $w$ by $\epsilon gt$ and then replacing $z$ by
$(g/a_2)u$ proves \eqref{eq:valuationtransform}.  For $r=1$,
$[u]\widehat F(u)^{-1}=-1$, giving $g^+=\epsilon a_2$.
\end{proof}

\section{A 2-adic nonvanishing theorem}\label{sec:valuation}

We regard $\mathbb Q$ as a subfield of the field $\Qtwo$ of 2-adic
numbers.  The valuation $v_2$ is normalized by $v_2(2)=1$ and is
extended by the convention $v_2(0)=+\infty$.  The ring of 2-adic
integers is denoted by $\Ztwo$.  For a rational number $x$,
\[
 x\in\Ztwo\Longleftrightarrow v_2(x)\geq0,
 \qquad
 x\in2\Ztwo\Longleftrightarrow v_2(x)\geq1.
\]
The passage from $\mathbb Q$ to $\Qtwo$ changes none of the rational
coefficients; it records their divisibility by powers of $2$.

\begin{theorem}[Nonvanishing at every depth]\label{thm:nonzero}
The series $A_i$ in \eqref{eq:iteration} are defined for every
$i\geq1$, and
\begin{equation}\label{eq:valuation}
  v_2(g_i)=1-2^{i-2}\qquad(i\geq2).
\end{equation}
In particular, every $g_i$ is a nonzero rational number.
\end{theorem}

\begin{proof}
Consider one stage
\begin{equation}\label{eq:oneStage}
  A(z)=gz+a_2z^2+\sum_{k\geq3}a_kz^k
\end{equation}
with $g\ne0$ and $a_2\ne0$, and put
\begin{equation}\label{eq:qdef}
  q=v_2(a_2/g).
\end{equation}
We use the following three-part induction invariant:
\begin{equation}\label{eq:invariant}
  v_2(a_2)<1,
  \qquad q<0,
  \qquad
  v_2(a_{r+1}/g)>rq\quad(r\geq2).
\end{equation}

With $u=(a_2/g)z$, as defined in \eqref{eq:udef}, we have
\begin{equation}\label{eq:FhatExpansion}
 \widehat F(u)=1+u+\sum_{r\geq2}d_ru^r,
 \qquad
 d_r=\frac{a_{r+1}/g}{(a_2/g)^r}.
\end{equation}
The last condition in \eqref{eq:invariant} gives
\[
 v_2(d_r)=v_2(a_{r+1}/g)-rq>0.
\]
The valuation of a nonzero rational number is integral, so
$v_2(d_r)\geq1$; the same membership statement is automatic if
$d_r=0$.  Hence
\begin{equation}\label{eq:drEven}
  d_r\in2\Ztwo\qquad(r\geq2).
\end{equation}

Let $A^+=\mathcal L_\epsilon(A)$ and write its coefficients as in
\eqref{eq:Aplusdef}.  Lemma~\ref{lem:transformcoeff} gives
$g^+=\epsilon a_2$.  Set
\[
  h_2=[u^2]\widehat F(u)^{-2}.
\]
Using \eqref{eq:FhatExpansion},
\[
  h_2=3-2d_2.
\]
Thus $h_2$ is an odd 2-adic integer.  Formula
\eqref{eq:valuationtransform} with $r=2$ gives
\begin{equation}\label{eq:a2plus}
  a_2^+=-\frac{a_2^2h_2}{2}.
\end{equation}
Consequently,
\begin{equation}\label{eq:qplus}
  q^+:=v_2(a_2^+/g^+)=v_2(a_2)-1<0,
\end{equation}
where the strict inequality uses the first condition in
\eqref{eq:invariant}.  The same condition persists because
\begin{equation}\label{eq:a2plusbound}
  v_2(a_2^+)=2v_2(a_2)-1<1.
\end{equation}

For $r\geq2$, define
\begin{equation}\label{eq:Cdef}
  C_{r+1}=[u^{r+1}]\widehat F(u)^{-(r+1)}.
\end{equation}
Since $\widehat F(u)\in1+u\Ztwo[[u]]$, every negative integer power of
$\widehat F$ belongs to $\Ztwo[[u]]$.  Thus $C_{r+1}\in\Ztwo$ and
$v_2(C_{r+1})\geq0$.  No oddness assertion for $C_{r+1}$ is required;
the only coefficient whose oddness was used is $h_2$.  Applying
\eqref{eq:valuationtransform} with index $r+1$ and dividing by
$g^+=\epsilon a_2$ yields
\begin{align}
 v_2(a_{r+1}^+/g^+)
 &=r v_2(a_2)-v_2(r+1)+v_2(C_{r+1})\notag\\
 &\geq r v_2(a_2)-v_2(r+1)\notag\\
 &>r\bigl(v_2(a_2)-1\bigr)=rq^+,
 \label{eq:higherPersist}
\end{align}
because $v_2(r+1)<r$ for every $r\geq2$.  Equations
\eqref{eq:qplus}, \eqref{eq:a2plusbound}, and
\eqref{eq:higherPersist} prove that the complete invariant persists.

It remains to verify the initial stage $A_1=B$.  Here
\[
  g=b_1=2,
  \qquad
  a_2=b_2=-7,
  \qquad
  v_2(a_2)=0<1,
  \qquad
  q=v_2(-7/2)=-1.
\]
For $r\geq2$, Lemma~\ref{lem:stable} gives
$(r+1)b_{r+1}\in2\mathbb Z$, and hence
\[
 v_2(b_{r+1}/2)
 \geq-v_2(r+1)>-r=rq.
\]
Thus \eqref{eq:invariant} holds initially, and induction proves that
every transformed series exists with nonzero linear coefficient.

Finally, $g_2=7$, so $v_2(g_2)=0$.  At stage $i$, the identities
$g_{i+1}=\epsilon_i a_{2,i}$ and \eqref{eq:a2plus} give
\[
 v_2(g_{i+2})=2v_2(g_{i+1})-1\qquad(i\geq1).
\]
Solving this recurrence proves \eqref{eq:valuation}.
\end{proof}

\begin{corollary}[Doubly exponential denominator]\label{cor:denominator}
For $i\geq3$, the exact power of $2$ dividing the denominator of $g_i$
in lowest terms is $2^{2^{i-2}-1}$.
\end{corollary}

\begin{proof}
This is an immediate reformulation of \eqref{eq:valuation}.
\end{proof}

The corollary concerns reduced denominators.  By itself it gives no
lower bound on the ordinary real magnitude $|g_i|$, because a rational
number with a large denominator may be arbitrarily close to zero.

\section{The closed nested-sum formula}\label{sec:closed}

For integers $x>y\geq1$, define
\begin{equation}\label{eq:weights}
 f_1(x,y)=\frac{(-1)^y2^{x-2y}}{(x-y-1)!},
 \qquad
 f_i(x,y)=\frac{g_i^{x-2y}}{(x-y-1)!}\quad(i\geq2).
\end{equation}
All constants $g_i$ in \eqref{eq:weights} are nonzero by
Theorem~\ref{thm:nonzero}, so negative exponents cause no ambiguity.

For a strictly decreasing chain of positive integers
\begin{equation}\label{eq:chain}
  v_0>v_1>\cdots>v_r\geq1,
\end{equation}
define its terminal weight by
\begin{equation}\label{eq:zweight}
 z_{r+1}(v_0,\ldots,v_r)
 =\frac{g_{r+1}^{v_r}}{v_r!}
  \prod_{k=0}^{r-1}f_{k+1}(v_k,v_{k+1}).
\end{equation}

\begin{theorem}[Closed nested representation]\label{thm:full}
For every $s\geq1$, $m,n>s+1$, and $N=mn$, the number of matchings of
size $s$ in $G_{m,n}$ is
\begin{align}
 a(N,s)={}&\frac{2^s}{s!}N^s\notag\\
 &+\sum_{r=1}^{s-1}
   \sum_{s=v_0>v_1>\cdots>v_r\geq1}
   N^{s-v_1}z_{r+1}(v_0,\ldots,v_r)
   \left(\prod_{k=0}^{r-2}
      v_k^{v_{k+1}-v_{k+2}-1}\right)
   v_{r-1}^{v_r-1}.
 \label{eq:full}
\end{align}
For $r=1$, the product indexed by $k$ is empty and is interpreted as
$1$.  Formula \eqref{eq:full} is finite and has maximum nesting depth
$s-1$.
\end{theorem}

\begin{proof}
By \eqref{eq:stableP},
\[
  a(N,s)=[z^s]e^{NB(z)}.
\]
The coefficient of $N^s$ is $g_1^s/s!=2^s/s!$.
Fix $v_1$ with $s>v_1\geq1$ and take the coefficient of
$N^{s-v_1}$.  Apply Lemma~\ref{lem:lagrange} with
\[
  A=A_1,
  \qquad n=s,
  \qquad j=v_1,
  \qquad \epsilon=\epsilon_1=-1.
\]
The prefactor is $f_1(s,v_1)$.  In the sum
\eqref{eq:coeffsecond}, the term with $\ell=0$ terminates the
transformation after the first level and contributes
\[
  f_1(s,v_1)s^{v_1-1}\frac{g_2^{v_1}}{v_1!}.
\]

For a term with $\ell=v_2\geq1$, where $v_1>v_2$, the remaining
coefficient is
\[
  [t^{v_1}]
  \frac{A_2(t)^{v_1-v_2}}{(v_1-v_2)!}
  =[X^{v_1-v_2}][t^{v_1}]e^{XA_2(t)}.
\]
Apply Lemma~\ref{lem:lagrange} again, now with
\[
  A=A_2,
  \qquad n=v_1,
  \qquad j=v_2,
  \qquad \epsilon=\epsilon_2=1.
\]
This produces the prefactor $f_2(v_1,v_2)$ and the power
$s^{v_1-v_2-1}$ carried by the preceding level.

Repeating the argument at
\[
  n=v_2,v_3,\ldots,v_{r-1}
\]
produces the successive factors
$f_3(v_2,v_3),\ldots,f_r(v_{r-1},v_r)$ and the powers displayed in
\eqref{eq:full}.  Every application strictly decreases a positive
integer index, so the process terminates after at most $s-1$ steps.
At the terminal index $v_r$, the linear term of $A_{r+1}$ contributes
\[
  \frac{g_{r+1}^{v_r}}{v_r!}
  v_{r-1}^{v_r-1}.
\]
Combining the prefactors gives precisely
$z_{r+1}(v_0,\ldots,v_r)$ with $v_0=s$, and summing over all chains
proves \eqref{eq:full}.
\end{proof}

\begin{corollary}[Number of structural contributions]\label{cor:terms}
Formula \eqref{eq:full} has $2^{s-1}-1$ nonleading chain
contributions and $2^{s-1}$ contributions in total.
\end{corollary}

\begin{proof}
A decreasing chain beginning at $s$ is uniquely determined by its
nonempty set of remaining entries, which is a subset of
$\{1,\ldots,s-1\}$.  There are $2^{s-1}-1$ such subsets.  Adding the
leading term proves the assertion.
\end{proof}

\begin{remark}[Extension beyond dimers]\label{rem:generality}
The nested construction is not intrinsically a dimer construction.
Suppose a periodic finite-range hard-core model has, in its stable
range, a connected series
\[
 B_{\mathcal M}(z)=\sum_{k\geq1}b_k^{\mathcal M}z^k
 \quad\hbox{and}\quad
 P_s^{\mathcal M}(N)=[z^s]\exp\!\bigl(NB_{\mathcal M}(z)\bigr).
\]
Apply the transformations of Section~\ref{sec:transform} to
$B_{\mathcal M}$.  If their successive linear coefficients
$g_i^{\mathcal M}$ are nonzero, the proof of
Theorem~\ref{thm:full} applies verbatim: in
\eqref{eq:weights}, replace $2=g_1$ by
$g_1^{\mathcal M}$ and every later $g_i$ by
$g_i^{\mathcal M}$.  This covers, in particular, the stable
fixed-population polynomials of aligned $k$-mers for arbitrary
$k\geq2$, as well as other finite-range hard-core placement models
covered by the polynomiality theorem of \cite{Kong2026}, subject only
to this explicit nondegeneracy condition.  Thus the formal nested
formula is universal; the 2-adic proof of nondegeneracy in
Theorem~\ref{thm:nonzero} is the model-specific part for square-lattice
dimers.  Establishing analogous arithmetic criteria for other models is
a natural next problem.
\end{remark}

For illustration, the first four stable polynomials obtained either
from \eqref{eq:Pdef} or from \eqref{eq:full} are
\begin{align*}
 P_1(N)&=2N,\\
 P_2(N)&=2N^2-7N,\\
 P_3(N)&=\frac43N^3-14N^2+\frac{116}{3}N,\\
 P_4(N)&=\frac23N^4-14N^3+\frac{611}{6}N^2-\frac{521}{2}N.
\end{align*}
The next constants generated by \eqref{eq:iteration} are
\begin{equation}\label{eq:firstgvalues}
  g_3=\frac{23}{6},
  \qquad
  g_4=\frac{2831}{24},
  \qquad
  g_5=-\frac{1219727477}{155520}.
\end{equation}

\section{Exact coefficients and numerical growth}\label{sec:numerics}

The density series associated with the stable Mayer series is
\begin{equation}\label{eq:density}
  zB'(z)=\sum_{k\geq1}c_kz^k,
  \qquad c_k=kb_k.
\end{equation}
Table~\ref{tab:ck} gives the exact square-lattice coefficients through
$k=26$.  They were obtained by independent exact computations using a
transfer-matrix method and a cavity-ratio method, and the outputs were
cross-checked against one another.  For $1\leq k\leq24$, the values
coincide term by term with the square-lattice column of Table~I in
Butera and Pernici \cite{ButeraPernici2012}.  In that paper the
coefficients $d_k$ are defined by the dimer density per site
$\rho(z)=\sum_{k\geq1}d_kz^k$.  Equation~\eqref{eq:density} is the same
density series, so the conventions agree exactly and $d_k=c_k$;
The lower-order terms also agree
indirectly with Gaunt's activity expansion \cite{Gaunt1969}.
The last two entries extend
the published density sequence.  Every entry is even, as proved in
Lemma~\ref{lem:stable}.

\begin{table}[ht]
\centering
\caption{Exact coefficients $c_k=kb_k$ through order $26$.}
\label{tab:ck}
\small
\begin{tabular}{r@{\quad}r@{\qquad}r@{\quad}r}
\toprule
$k$ & $c_k$ & $k$ & $c_k$\\
\midrule
1  & $2$                         & 14 & $-12164408791920$\\
2  & $-14$                       & 15 & $129177146454536$\\
3  & $116$                       & 16 & $-1376741271026898$\\
4  & $-1042$                     & 17 & $14719835348283940$\\
5  & $9812$                      & 18 & $-157827022198103624$\\
6  & $-95288$                    & 19 & $1696509872615736256$\\
7  & $945688$                    & 20 & $-18277500804075889672$\\
8  & $-9537906$                  & 21 & $197318358676755340504$\\
9  & $97398764$                  & 22 & $-2134157389758234716560$\\
10 & $-1004479624$               & 23 & $23121845343883130936248$\\
11 & $10442811216$               & 24 & $-250895842743288634987656$\\
12 & $-109291830952$             & 25 & $2726369642483654704568312$\\
13 & $1150263509280$             & 26 & $-29665252708759199866230352$\\
\bottomrule
\end{tabular}
\end{table}

The transformation \eqref{eq:iteration} was evaluated through $i=26$
with exact arbitrary-precision rational arithmetic, using multiple
independent implementations.  Table~\ref{tab:growth} records
diagnostics for $g_i$.  The columns $D_i^{\rm num}$ and $D_i^{\rm den}$
are the decimal digit counts of the reduced numerator and denominator.
The last two columns are
\[
  L_i=\log\log|g_i|,
  \qquad
  R_i=|g_i|^{1/2^i}.
\]
The decimals are derived only after the exact rational computation.

\begingroup
\scriptsize
\setlength{\tabcolsep}{5pt}
\begin{longtable}{r r r r r}
\caption{Growth diagnostics for the exact rational constants $g_i$
through $i=26$.}\label{tab:growth}\\
\toprule
$i$ & $D_i^{\rm num}$ & $D_i^{\rm den}$ & $L_i$ & $R_i$\\
\midrule
1  & 1     & 1     & $-0.3665$ & 1.41421\\
2  & 1     & 1     & 0.6657  & 1.62658\\
3  & 2     & 1     & 0.2955  & 1.18290\\
4  & 4     & 2     & 1.5624  & 1.34736\\
5  & 10    & 6     & 2.1936  & 1.32343\\
6  & 18    & 11    & 2.8567  & 1.31251\\
7  & 38    & 23    & 3.5819  & 1.32419\\
8  & 75    & 45    & 4.2564  & 1.31734\\
9  & 153   & 92    & 4.9437  & 1.31520\\
10 & 307   & 183   & 5.6506  & 1.32021\\
11 & 616   & 369   & 6.3445  & 1.32050\\
12 & 1230  & 737   & 7.0339  & 1.31913\\
13 & 2464  & 1476  & 7.7287  & 1.31974\\
14 & 4925  & 2951  & 8.4221  & 1.31980\\
15 & 9855  & 5907  & 9.1150  & 1.31973\\
16 & 19709 & 11813 & 9.8082  & 1.31974\\
17 & 39421 & 23628 & 10.5013 & 1.31973\\
18 & 78841 & 47256 & 11.1944 & 1.31973\\
19 & 157682 & 94514 & 11.8876 & 1.31973\\
20 & 315361 & 189024 & 12.5807 & 1.31973\\
21 & 630732 & 378057 & 13.2739 & 1.31973\\
22 & 1261464 & 756114 & 13.9670 & 1.31973\\
23 & 2522929 & 1512230 & 14.6602 & 1.31973\\
24 & 5045857 & 3024459 & 15.3533 & 1.31973\\
25 & 10091716 & 6048920 & 16.0465 & 1.31973\\
26 & 20183432 & 12097839 & 16.7396 & 1.31973\\
\bottomrule
\end{longtable}
\endgroup

A least-squares fit over $20\leq i\leq26$ gives
\begin{equation}\label{eq:fit}
  \log\log|g_i|
  \approx0.693148\,i-1.28221,
\end{equation}
with coefficient of determination exceeding $0.999999999999$.  The fitted
slope is close to $\log2=0.693147\ldots$.  Equivalently, the data are
consistent with
\[
  \log|g_i|\sim\gamma 2^i,
  \qquad
  |g_i|^{1/2^i}\longrightarrow e^\gamma,
\]
where the last four values of $R_i$ agree to six decimal places and
give the observed limiting scale $1.3197289$.
The range is necessarily finite, and neither the quality of the fit nor
Corollary~\ref{cor:denominator} proves an Archimedean lower bound.
We therefore state the observation as a conjecture.

\begin{conjecture}[Doubly exponential magnitude]\label{conj:growth}
There exists a constant $\Lambda>1$ such that
\[
  \lim_{i\to\infty}|g_i|^{1/2^i}=\Lambda.
\]
The current exact data suggest $\Lambda\approx1.3197289$.
\end{conjecture}

\section{Complexity perspective and open problems}\label{sec:complexity}

Corollary~\ref{cor:terms} gives an exact combinatorial count of the
structural contributions in the closed formula.  If
$g_1,\ldots,g_s$ are already available, a straightforward evaluation
uses exponentially many rational terms.  The total length of all chain
products is of order $s2^s$.  This is an arithmetic-operation count
rather than a bit-complexity theorem: the numerators and denominators in
Table~\ref{tab:growth} show that operand sizes themselves can grow very
rapidly.

This distinction is also important when comparing the formula with
\#P-completeness.  Jerrum's theorem concerns the uniform problem of
counting matchings in planar graphs \cite{Jerrum1987}.  Here the graph
family is highly structured, the formula applies in the stable regime,
and the number $s$ of dimers is the parameter controlling the nesting
depth.  The present result neither separates FP from \#P nor supplies a
complete bit-complexity classification.  It does, however, expose an
exact algebraic structure on which such an analysis can be based.

Several questions arise naturally.
\begin{enumerate}
  \item Prove or disprove Conjecture~\ref{conj:growth}, and determine
  the constant $\Lambda$ if the limit exists.
  \item Determine whether the nested sum admits cancellations or a
  dynamic evaluation that reduces its arithmetic complexity.
  \item Quantify the bit complexity of producing $g_s$ from
  $c_1,\ldots,c_s$ and of evaluating \eqref{eq:full} exactly.
  \item In the general framework of Remark~\ref{rem:generality}, find
  model-specific criteria guaranteeing nonvanishing of every
  $g_i^{\mathcal M}$, particularly for aligned $k$-mers and other
  periodic finite-range hard-core models.
  \item Find a direct combinatorial interpretation of the constants
  $g_i$ or of their numerators and denominators.
\end{enumerate}

\section{Conclusion}\label{sec:conclusion}

The fixed-dimer canonical partition function on a sufficiently large
square torus is not merely a polynomial in the number of sites.  Its
coefficients possess a recursively nested structure that can be unfolded
completely by formal reversion and Lagrange inversion.  The resulting
finite formula is indexed by decreasing integer chains and is valid for
arbitrary $s$.

The reversion--logarithm mechanism is algebraic rather than
dimer-specific.  Once a finite-range hard-core model has a stable Mayer
series and a nondegenerate transformation sequence, the same decreasing
chain construction yields its fixed-population polynomials.  The square
monomer--dimer model supplies the first case in which nondegeneracy is
proved at every depth by an explicit arithmetic invariant.

The proof also reveals an arithmetic feature invisible in the original
matching counts.  The connected Mayer coefficients satisfy
$kb_k\in2\mathbb Z$, and this divisibility forces an exact doubly
exponential power of $2$ in the denominators of the transformation
constants.  The much stronger numerical phenomenon suggested by the
real magnitudes of these constants remains open.  Together, the closed
formula, the valuation theorem, and the exact data provide a rigorous
foundation for future work on the arithmetic and computational
complexity of monomer--dimer enumeration.

\end{document}